\documentclass[12pt]{article}                                               

\input{amssym.tex}                                                           
                                                                             
\begin{document}                                                             
\title{On algebraic values of function  $\exp ~(2\pi i ~x+\log\log y)$}

\author{Igor  ~Nikolaev}


\date{}
 \maketitle


\newtheorem{thm}{Theorem}
\newtheorem{lem}{Lemma}
\newtheorem{dfn}{Definition}
\newtheorem{rmk}{Remark}
\newtheorem{cor}{Corollary}
\newtheorem{prp}{Proposition}
\newtheorem{exm}{Example}
\newtheorem{cnj}{Conjecture}

\newcommand{\Qcoh}{\hbox{\bf Qcoh}}
\newcommand{\jac}{\hbox{\bf Jac}}
\newcommand{\ka}{\hbox{\bf k}}
\newcommand{\n}{\hbox{\bf n}}
\newcommand{\QGr}{\hbox{\bf QGr}}
\newcommand{\Z}{\hbox{\bf Z}}
\newcommand{\Q}{\hbox{\bf Q}}
\newcommand{\tr}{\hbox{\bf Trace}}

\begin{abstract}
It is  proved  that   for all but a finite set of the square-free integers $d$ 
 the value of transcendental  function $\exp~(2\pi i ~x+\log\log y)$
 is an  algebraic number for the algebraic  arguments  $x$ and $y$ lying in a 
 real  quadratic  field  of  discriminant  $d$.   Such a value generates the Hilbert class field 
of the imaginary quadratic field of  discriminant $-d$.

\vspace{7mm}

{\it Key words and phrases:  real  multiplication;  Sklyanin algebra}

\vspace{5mm}
{\it MSC: 11J81 (transcendence theory);   46L85 (noncommutative topology)}
\end{abstract}

\section{Introduction}
It is an old  problem
to determine  if  given irrational value of  a  transcendental function is 
algebraic or transcendental for certain algebraic arguments;  the algebraic values 
are particularly remarkable and worthy of thorough investigation,  see   [Hilbert  1902]  \cite{Hil1},  p. 456.
Only few general results are known,   see e.g.  [Baker  1975]  \cite{B}.   We shall mention the famous 
Gelfond-Schneider Theorem saying  that $e^{\beta\log\alpha}$ is a transcendental
number,  whenever $\alpha\not\in \{0, 1\}$ is an algebraic and $\beta$ an irrational
algebraic number.    In  contrast,   Klein's  invariant  $j(\tau)$   is known to take
algebraic values  whenever  $\tau\in {\Bbb H}:=\{x+iy\in {\Bbb C}~|~y>0\}$ is 
an imaginary quadratic number.

The aim of our  note is a result on  the algebraic  values  of   transcendental  function
\begin{equation}
{\cal J}(x,y):=\{e^{2\pi i ~x ~+ ~\log\log y} ~|-\infty<x<\infty, ~1<y<\infty\}
\end{equation}
for  the  arguments  $x$ and $y$ in a real quadratic field;  the function  ${\cal J}(x, y)$ 
can be viewed as an analog of  Klein's invariant $j(\tau)$,   hence  the  notation. 
Namely,  let ${\goth k}={\Bbb Q}(\sqrt{d})$ be a real quadratic field and ${\goth R}_{\goth f}=
{\Bbb Z}+{\goth f}O_{\goth k}$ be an order of conductor ${\goth f}\ge 1$ in the field ${\goth k}$;
let $h=|Cl~({\goth R}_{\goth f})|$ be the class number of ${\goth R}_{\goth f}$ and denote by
$\{{\Bbb Z}+{\Bbb Z}\theta_i ~|~ 1\le i\le h\}$ the set of pairwise non-isomorphic pseudo-lattices
in ${\goth k}$ having the same endomorphism ring ${\goth R}_{\goth f}$,  see [Manin 2004]  \cite{Man1},
Lemma 1.1.1. 
Finally, let $\varepsilon$ be the fundamental unit of ${\goth R}_{\goth f}$;  let $f\ge 1$ be the 
least integer satisfying equation  $|Cl~(R_f)|=|Cl~({\goth R}_{\goth f})|$,  where $R_f={\Bbb Z}+fO_k$
is an order of conductor $f$ in the imaginary quadratic field $k={\Bbb Q}(\sqrt{-d})$,  see \cite{Nik1}. 
Our main result can be formulated  as follows.   
\begin{thm}\label{thm1}
For each square-free positive integer  $d\not\in\{1,2,3,7,11,19,$
\linebreak
$43,67,163\}$ 
the values $\{{\cal J}(\theta_i,\varepsilon) ~|~ 1\le i\le h\}$  of transcendental 
function ${\cal J}(x,y)$  are algebraically conjugate  numbers  generating the 
Hilbert class field $H(k)$ of the imaginary quadratic field $k={\Bbb Q}(\sqrt{-d})$
modulo conductor $f$. 
\end{thm}
\begin{rmk}\label{rk1}
\textnormal{
Since  $H(k)\cong k(j(\tau))\cong {\Bbb Q}(f\sqrt{-d}, j(\tau))$ with $\tau\in R_f$,  one gets  an  
inclusion  ${\cal J}(\theta_i,\varepsilon))\in  {\Bbb Q}(f\sqrt{-d}, j(\tau))$. 
}
\end{rmk}
\begin{rmk}
\textnormal{
The absolute value $|z|=(z\bar z)^{1\over 2}$ of an algebraic
number $z$ is  always an ``abstract''  algebraic number,  i.e. a root of the
polynomial with integer coefficients;  yet   theorem \ref{thm1} implies  that
 $|{\cal J}(\theta_i,\varepsilon)|=\log\varepsilon$  is  a  transcendental number.   
This  apparent contradiction is false,   since 
quadratic extensions  of the field ${\Bbb Q}(z\bar z)$  have
no real embeddings in general;   in other words,   our  extension cannot be a subfield  
of  ${\Bbb R}$. 
}
\end{rmk}
The structure of article is as follows.  All preliminary facts
can be found   in  Section 2.    Theorem \ref{thm1} is proved in Section 3.

\section{Preliminaries}
The reader can find basics of the $C^*$-algebras in [Murphy 1990]  \cite{M}
and their $K$-theory  in [Blackadar  1986]  \cite{BL}.    
The noncommutative tori are covered in [Rieffel  1990]  \cite{Rie1}
and real multiplication in [Manin 2004]  \cite{Man1}.  
For  main ideas  of non-commutative algebraic geometry,  see the survey 
by [Stafford \& van ~den ~Bergh  2001]  \cite{StaVdb1}.

\subsection{Noncommutative tori}
By a  {\it noncommutative torus}  ${\cal A}_{\theta}$ one understands  the universal  {\it $C^*$-algebra}
generated by the unitary operators $u$ and $v$ acting on a Hilbert space ${\cal H}$
and satisfying the commutation relation $vu=e^{2\pi i\theta}uv$,  where $\theta$ is a 
real number.   
\begin{rmk}\label{rmk1}
\textnormal{
Note that  ${\cal A}_{\theta}$ is isomorphic to  a  free ${\Bbb C}$-algebra  on  four 
generators $u,u^*,v,v^*$   and  six   quadratic relations:
\begin{equation}\label{eq2}
\left\{
\begin{array}{cc}
vu &= e^{2\pi i\theta} uv,\\
v^*u^* &= e^{2\pi i\theta}u^*v^*,\\
v^*u &=  e^{-2\pi i\theta}uv^*,\\
vu^* &=  e^{-2\pi i\theta}u^*v,\\
u^*u &= uu^*=e,\\
v^*v &= vv^*=e. 
\end{array}
\right.
\end{equation}
Indeed,  the first and the last two relations in system (\ref{eq2})  are obvious 
from the definition of ${\cal A}_{\theta}$.   By way of example,   
let us demonstrate  that  relations $vu=e^{2\pi i\theta} uv$ and $u^*u=uu^*=v^*v=vv^*=e$ 
imply the relation  $v^*u =  e^{-2\pi i\theta}uv^*$ in system (\ref{eq2}).    Indeed, 
it follows from $uu^*=e$ and $vv^*=e$ that $uu^*vv^*=e$.  Since $uu^*=u^*u$ we can bring  the last 
equation to  the form  $u^*uvv^*=e$ and multiply the  both sides by the constant $e^{2\pi i\theta}$;
thus one gets the equation $u^*(e^{2\pi i\theta}uv)v^*=e^{2\pi i\theta}$.  
But $e^{2\pi i\theta}uv=vu$ and our main equation takes the form $u^*vuv^*= e^{2\pi i\theta}$. 
We can multiply on the left the  both sides of the equation by the element $u$
and thus get the equation $uu^*vuv^*= e^{2\pi i\theta}u$; since $uu^*=e$ 
one  arrives at the equation  $vuv^*= e^{2\pi i\theta}u$.  
Again one can multiply on the left the both sides  by the element 
$v^*$ and thus get the equation  $v^*vuv^*= e^{2\pi i\theta}v^*u$;   since $v^*v=e$
one gets  $uv^*= e^{2\pi i\theta}v^*u$ and the required identity $v^*u =  e^{-2\pi i\theta}uv^*$.
The remaining two relations in (\ref{eq2})  are proved likewise;   we leave it   to the reader as an exercise in
non-commutative algebra.   
}
\end{rmk}

\bigskip
Recall that  the algebra ${\cal A}_{\theta}$ is said to be {\it stably isomorphic}
(Morita equivalent)  to ${\cal A}_{\theta'}$,   whenever ${\cal A}_{\theta}\otimes {\cal K}\cong 
{\cal A}_{\theta'}\otimes {\cal K}$,  where ${\cal K}$ is the $C^*$-algebra of all compact operators
on ${\cal H}$;  the ${\cal A}_{\theta}$ is stably isomorphic to ${\cal A}_{\theta'}$ if and only if 
\begin{equation}\label{eq3}
\theta'={a\theta +b\over c\theta+d}\quad
\hbox{for some matrix} \quad  \left(\matrix{a & b\cr c & d}\right)\in SL_2({\Bbb Z}). 
\end{equation}
The $K$-theory of ${\cal A}_{\theta}$ is  two-periodic and 
$K_0({\cal A}_{\theta})\cong  K_1({\cal A}_{\theta})\cong {\Bbb Z}^2$ so that  
the Grothendieck  semigroup $K_0^+({\cal A}_{\theta})$ corresponds  to positive reals of 
the set   ${\Bbb Z}+{\Bbb Z}\theta\subset {\Bbb R}$ called a {\it pseudo-lattice}.  
The torus ${\cal A}_{\theta}$ is said to have {\it real multiplication},    if $\theta$ is a quadratic
irrationality,  i.e.  irrational root of a quadratic polynomial with integer coefficients.  
The real multiplication says that the endomorphism ring of pseudo-lattice 
${\Bbb Z}+{\Bbb Z}\theta$ exceeds the ring ${\Bbb Z}$ corresponding to multiplication
by $m$ endomorphisms;  similar to complex multiplication, it means that the
endomorphism ring is isomorphic to an order ${\goth R}_{\goth f}={\Bbb Z}+{\goth f}O_{\goth k}$
of conductor ${\goth f}\ge 1$ in the real quadratic field ${\goth k}={\Bbb Q}(\theta)$ -- 
hence the name,   see [Manin 2004]  \cite{Man1}.   If $d>0$ is the discriminant of  ${\goth k}$,   then 
by ${\cal A}_{RM}^{(d, {\goth f})}$   we denote a noncommutative torus  with real multiplication
by the order ${\goth R}_{\goth f}$.

\subsection{Elliptic curves}
For the sake of clarity,  let us recall some well-known facts.   An {\it elliptic curve}
is the subset of the complex projective plane of the form
${\cal E}({\Bbb C})=\{(x,y,z)\in {\Bbb C}P^2 ~|~ y^2z=4x^3+axz^2+bz^3\}$,
where $a$ and  $b$  are some constant complex numbers.   Recall that  one can embed 
${\cal E}({\Bbb C})$ into the complex projective space ${\Bbb C}P^3$  as the set
of points of intersection of two {\it quadric surfaces}  given by the  system of homogeneous 
equations 
\begin{equation}\label{eq4}
\left\{
\begin{array}{ccc}
u^2+v^2+w^2+z^2 &=&  0,\\
Av^2+Bw^2+z^2  &=&  0,   
\end{array}
\right.
\end{equation}
where $A$ and $B$ are some constant  complex numbers and  
$(u,v,w,z)\in {\Bbb C}P^3$;  the system (\ref{eq4})  is called
the {\it Jacobi form} of elliptic curve ${\cal E}({\Bbb C})$.   
  Denote by  ${\Bbb H}=\{x+iy\in {\Bbb C}~|~y>0\}$  the Lobachevsky
 half-plane;  whenever  $\tau\in {\Bbb H}$,  one gets a complex torus  ${\Bbb C}/({\Bbb Z}+{\Bbb Z}\tau)$.
 Each complex torus is isomorphic  to  a non-singular elliptic curve;  the isomorphism  is realized by the
Weierstrass $\wp$ function and we shall write ${\cal E}_{\tau}$ to denote the corresponding elliptic curve.
Two elliptic curves ${\cal E}_{\tau}$ and ${\cal E}_{\tau'}$ are isomorphic if and only if 
\begin{equation}\label{eq5}
\tau'={a\tau +b\over c\tau+d}\quad
\hbox{for some matrix} \quad  \left(\matrix{a & b\cr c & d}\right)\in SL_2({\Bbb Z}). 
\end{equation}
If $\tau$ is an imaginary quadratic number,  elliptic curve ${\cal E}_{\tau}$  is said to have 
{\it complex multiplication};   in this case lattice ${\Bbb Z}+{\Bbb Z}\tau$
admits non-trivial endomorphisms realized as multiplication of points of the lattice by the  imaginary
quadratic numbers,  hence the name.  We shall write ${\cal E}_{CM}^{(-d,f)}$ to denote elliptic curve with complex
multiplication by an order $R_f={\Bbb Z}+fO_k$  of conductor $f\ge 1$ in the imaginary quadratic
field $k={\Bbb Q}(\sqrt{-d})$.

\subsection{Sklyanin algebras}
By the  {\it Sklyanin algebra} $S_{\alpha,\beta,\gamma}({\Bbb C})$  one understands  a 
free   ${\Bbb C}$-algebra   on   four  generators    and  six   relations: 
\begin{equation}
\left\{
\begin{array}{ccc}
x_1x_2-x_2x_1 &=& \alpha(x_3x_4+x_4x_3),\\
x_1x_2+x_2x_1 &=& x_3x_4-x_4x_3,\\
x_1x_3-x_3x_1 &=& \beta(x_4x_2+x_2x_4),\\
x_1x_3+x_3x_1 &=& x_4x_2-x_2x_4,\\
x_1x_4-x_4x_1 &=& \gamma(x_2x_3+x_3x_2),\\ 
x_1x_4+x_4x_1 &=& x_2x_3-x_3x_2,
\end{array}
\right.
\end{equation}
where $\alpha+\beta+\gamma+\alpha\beta\gamma=0$.  The algebra  $S_{\alpha,\beta,\gamma}({\Bbb C})$
represents   a   twisted homogeneous {\it  coordinate ring}   of an elliptic curve  ${\cal E}_{\alpha,\beta,\gamma}({\Bbb C})$
given in  its  Jacobi form
\begin{equation}
\left\{
\begin{array}{ccc}
u^2+v^2+w^2+z^2 &=&  0,\\
{1-\alpha\over 1+\beta}v^2+
 {1+\alpha\over 1-\gamma}w^2+z^2  &=&  0,   
\end{array}
\right.
\end{equation}
see  [Smith \& Stafford 1993]  \cite{SmiSta1}, p.267  and   
[Stafford \& van ~den ~Bergh  2001]  \cite{StaVdb1},   Example 8.5.
The latter  means that algebra $S_{\alpha,\beta,\gamma}({\Bbb C})$ satisfies an isomorphism
$\hbox{{\bf Mod}}~(S_{\alpha,\beta,\gamma}({\Bbb C}))/
\hbox{{\bf Tors}}\cong \hbox{{\bf Coh}}~({\cal E}_{\alpha,\beta,\gamma}({\Bbb C}))$,
 where {\bf Coh} is  the category of quasi-coherent sheaves on ${\cal E}_{\alpha,\beta,\gamma}({\Bbb C})$, 
  {\bf Mod}  the category of graded left modules over the graded ring $S_{\alpha,\beta,\gamma}({\Bbb C})$
 and  {\bf Tors}  the full sub-category of {\bf Mod} consisting of the
torsion modules,  see  [Stafford \& van ~den ~Bergh  2001]  \cite{StaVdb1},  p.173.  
The algebra $S_{\alpha,\beta,\gamma}({\Bbb C})$ defines  a natural {\it automorphism} 
$\sigma$ of  elliptic curve ${\cal E}_{\alpha,\beta,\gamma}({\Bbb C})$,   {\it ibid.}

\section{Proof of theorem \ref{thm1}}
For the sake of clarity,  let us outline  main ideas.    The proof is based on a 
categorical correspondence (a covariant functor) between elliptic curves ${\cal E}_{\tau}$ 
and noncommutative tori ${\cal A}_{\theta}$ taken with their  ``scaled units''
${1\over\mu}e$.  Namely,  we prove that for $\sigma^4=Id$  the norm-closure of a self-adjoint 
representation of the Sklyanin algebra $S_{\alpha,\beta,\gamma}({\Bbb C})$ 
by the linear operators $u=x_1,u^*=x_2, v=x_3, v^*=x_4$ on a Hilbert space 
${\cal H}$ is isomorphic to the $C^*$-algebra ${\cal A}_{\theta}$ so that its 
unit $e$ is scaled by a positive real $\mu$, see lemma \ref{lem2};   because  $S_{\alpha,\beta,\gamma}({\Bbb C})$
is a coordinate ring of elliptic curve ${\cal E}_{\alpha,\beta,\gamma}({\Bbb C})$  
so will be the algebra ${\cal A}_{\theta}$   modulo the unit ${1\over\mu}e$.  
Moreover,  our construction  entails  that a  coefficient  $q$ of  elliptic curve 
 ${\cal E}_{\alpha,\beta,\gamma}({\Bbb C})$ is linked to  the  constants 
$\theta$ and $\mu$   by the formula $q=\mu e^{2\pi i\theta}$, see lemma \ref{lem1}.    
Suppose that our elliptic curve has complex multiplication, i.e.  ${\cal E}_{\tau}\cong {\cal E}_{CM}^{(-d,f)}$;
then its coordinate ring $({\cal A}_{\theta},  {1\over\mu}e)$  must have real multiplication,  i.e.  ${\cal A}_{\theta}\cong {\cal A}_{RM}^{(d, {\goth f})}$
and $\mu=\log\varepsilon$, where  $|Cl~(R_f)|=|Cl~({\goth R}_{\goth f})|$ and $\varepsilon$ is the fundamental unit of 
order ${\goth R}_{\goth f}$,   see lemma \ref{lem3}.   But elliptic curve  ${\cal E}_{CM}^{(-d,f)}$
has coefficients in  the Hilbert class field $H(k)$ over imaginary quadratic field $k={\Bbb Q}(\sqrt{-d})$ modulo conductor $f$;
thus $q\in H(k)$ and  therefore one gets an inclusion
\begin{equation}
\mu e^{2\pi i\theta}\in H(k),
\end{equation}
where $\theta\in {\goth k}= {\Bbb Q}(\sqrt{d})$ and $\mu=\log\varepsilon$.  
(Of course,  our argument  is valid only when $q\not\in {\Bbb R}$,  i.e. when  $|Cl~(R_f)|\ge 2$;
but there are only a finite number of discriminants $d$ with  $|Cl~(R_f)|=1$.)   
Let us pass  to a detailed argument.
\begin{lem}\label{lem1}
If $\sigma^4=Id$,  then the Sklyanin algebra  $S_{\alpha,\beta,\gamma}({\Bbb C})$
endowed with the involution $x_1^*=x_2$ and $x_3^*=x_4$ 
is isomorphic to  a free algebra ${\Bbb C}\langle x_1,x_2,x_3,x_4\rangle$  modulo an ideal 
generated  by  six   quadratic  relations
\begin{equation}\label{eq9}
\left\{
\begin{array}{cc}
x_3x_1 &= \mu e^{2\pi i\theta}x_1x_3,\\
x_4x_2 &= {1\over \mu} e^{2\pi i\theta}x_2x_4,\\
x_4x_1 &= \mu e^{-2\pi i\theta}x_1x_4,\\
x_3x_2 &= {1\over \mu} e^{-2\pi i\theta}x_2x_3,\\
x_2x_1 &= x_1x_2,\\
x_4x_3 &= x_3x_4,
\end{array}
\right.
\end{equation}
where $\theta=Arg~(q)$ and $\mu=|q|$ for a complex number $q\in {\Bbb C}\setminus \{0\}$. 
\end{lem}
{\it Proof.}  
(i)  Since  $\sigma^4=Id$,   the Sklyanin algebra $S_{\alpha, \beta, \gamma}({\Bbb C})$
is isomorphic to a free   algebra  ${\Bbb C}\langle x_1,x_2,x_3,x_4\rangle$
modulo an ideal generated by   the  skew-symmetric relations 
\begin{equation}\label{eq10}
\left\{
\begin{array}{ccc}
x_3x_1 &=& q_{13} x_1x_3,\\
x_4x_2 &=&  q_{24}x_2x_4,\\
x_4x_1 &=&  q_{14}x_1x_4,\\
x_3x_2 &=&  q_{23}x_2x_3,\\
x_2x_1&=&  q_{12}x_1x_2,\\
x_4x_3&=&  q_{34}x_3x_4,
\end{array}
\right.
\end{equation}
where $q_{ij}\in {\Bbb C}\setminus\{0\}$,  see  [Feigin \& Odesskii  1989]  \cite{FeOd1},  Remark 1.

\bigskip
(ii) It is verified directly,  that  relations (\ref{eq10})   are invariant of the involution 
$x_1^*=x_2$ and  $x_3^*=x_4$,     if and only if
\begin{equation}\label{eq11}
\left\{
\begin{array}{ccc}
q_{13} &=&  (\bar q_{24})^{-1},\\
q_{24} &=&  (\bar q_{13})^{-1},\\
q_{14} &= & (\bar q_{23})^{-1},\\
q_{23} &= & (\bar q_{14})^{-1},\\
q_{12} &= & \bar q_{12},\\
q_{34} &= & \bar q_{34},
\end{array}
\right.
\end{equation}
where $\bar q_{ij}$ means the complex conjugate of $q_{ij}\in {\Bbb C}\setminus\{0\}$.

\bigskip
\begin{rmk}
\textnormal{
The invariant  relations (\ref{eq11})   define an involution on the Sklyanin algebra;  
we shall refer to such as  a {\it Sklyanin $\ast$-algebra}.
}
\end{rmk}

\bigskip
(iii)
Consider a one-parameter family  $S(q_{13})$ of the Sklyanin $\ast$-algebras
defined by the following additional constraints
\begin{equation}
\left\{
\begin{array}{ccc}
q_{13} &=& \bar  q_{14},\\
q_{12} &=&  q_{34}=1.
\end{array}
\right.
\end{equation}
It is not hard to see,  that the $\ast$-algebras  $S(q_{13})$ 
are pairwise non-isomorphic for different values of complex  parameter $q_{13}$;
therefore  the family $S(q_{13})$ is   a  normal form of  the Sklyanin $\ast$-algebra 
$S_{\alpha, \beta, \gamma}({\Bbb C})$ with $\sigma^4=Id$. 
It remains to notice, that one can write  complex parameter $q:=q_{13}$
in the polar form $q=\mu e^{2\pi i\theta}$,  where $\theta=Arg~(q)$
and $\mu=|q|$.   Lemma \ref{lem1} follows.
$\square$

\begin{lem}\label{lem2}
{\bf (basic isomorphism)}
The system of relations (\ref{eq2})  for noncommutative torus ${\cal A}_{\theta}$
with $u=x_1, u^*=x_2, v=x_3, v^*=x_4$,  i.e.  
\begin{equation}\label{eq13}
\left\{
\begin{array}{cc}
x_3x_1 &=  e^{2\pi i\theta}x_1x_3,\\
x_4x_2 &=  e^{2\pi i\theta}x_2x_4,\\
x_4x_1 &=  e^{-2\pi i\theta}x_1x_4,\\
x_3x_2 &=   e^{-2\pi i\theta}x_2x_3,\\
x_2x_1 &= x_1x_2=e,\\
x_4x_3 &= x_3x_4=e, 
\end{array}
\right.
\end{equation}
is equivalent to the system of relations (\ref{eq9}) for the Sklyanin $\ast$-algebra, i.e.
\begin{equation}\label{eq14}
\left\{
\begin{array}{cc}
x_3x_1 &= \mu e^{2\pi i\theta}x_1x_3,\\
x_4x_2 &= {1\over \mu} e^{2\pi i\theta}x_2x_4,\\
x_4x_1 &= \mu e^{-2\pi i\theta}x_1x_4,\\
x_3x_2 &= {1\over \mu} e^{-2\pi i\theta}x_2x_3,\\
x_2x_1 &= x_1x_2,\\
x_4x_3 &= x_3x_4,
\end{array}
\right.
\end{equation}
modulo the following  ``scaled unit relation'' 
\begin{equation}\label{eq15}
x_1x_2=x_3x_4={1\over\mu}e.
\end{equation}
\end{lem}
{\it Proof.} 
(i) Using the last two relations,  one can bring the noncommutative torus 
relations (\ref{eq13})   to the form
\begin{equation}\label{eq16}
\left\{
\begin{array}{ccc}
x_3x_1x_4 &=&  e^{2\pi i\theta}x_1,\\
x_4 &= & e^{2\pi i\theta}x_2x_4x_1,\\
x_4x_1x_3 &=&  e^{-2\pi i\theta}x_1,\\
x_2 &=&  e^{-2\pi i\theta}x_4x_2x_3,\\
x_1x_2   &=&   x_2x_1   =e,\\
 x_3x_4    &=&  x_4x_3  =e.
\end{array}
\right.
\end{equation}

\bigskip
(ii)  The system of relations (\ref{eq14})  for the Sklyanin $\ast$-algebra complemented 
by the scaled unit relation (\ref{eq15}),  i.e.  
\begin{equation}\label{eq17}
\left\{
\begin{array}{cc}
x_3x_1 &= \mu e^{2\pi i\theta}x_1x_3,\\
x_4x_2 &= {1\over \mu} e^{2\pi i\theta}x_2x_4,\\
x_4x_1 &= \mu e^{-2\pi i\theta}x_1x_4,\\
x_3x_2 &= {1\over \mu} e^{-2\pi i\theta}x_2x_3,\\
x_2x_1 &= x_1x_2={1\over\mu}e,\\
x_4x_3 &= x_3x_4={1\over\mu}e,
\end{array}
\right.
\end{equation}
is equivalent to the system
\begin{equation}\label{eq18}
\left\{
\begin{array}{cc}
x_3x_1x_4 &= e^{2\pi i\theta}x_1,\\
x_4 &= e^{2\pi i\theta}x_2x_4x_1,\\
x_4x_1x_3 &= e^{-2\pi i\theta}x_1,\\
x_2 &= e^{-2\pi i\theta}x_4x_2x_3,\\
x_2x_1 &= x_1x_2={1\over\mu}e,\\
x_4x_3 &= x_3x_4={1\over\mu}e
\end{array}
\right.
\end{equation}
by using  multiplication and cancellation involving the last two equations.

\bigskip
(iii)  For each $\mu\in (0,\infty)$ consider a {\it scaled unit}  $e':={1\over\mu} e$ of
the Sklyanin $\ast$-algebra $S(q)$ and the two-sided ideal $I_{\mu}\subset S(q)$
generated by the relations $x_1x_2=x_3x_4=e'$.    Comparing the defining relations (\ref{eq14}) for
$S(q)$ with relation (\ref{eq13}) for the noncommutative torus ${\cal A}_{\theta}$,   one gets an
isomorphism 
\begin{equation}\label{eq19}
S(q)~/~I_{\mu}\cong {\cal A}_{\theta}. 
\end{equation}
The isomorphism maps  generators $x_1,\dots,x_4$ of  the Sklyanin
$\ast$-algebra $S(q)$ to such  of the $C^*$-algebra ${\cal A}_{\theta}$ and 
the {\it scaled} unit $e'\in S(q)$ to the {\it ordinary} unit $e\in {\cal A}_{\theta}$. 
Lemma \ref{lem2} follows.  
$\square$

\begin{rmk}
\textnormal{
It follows from (\ref{eq19})  that noncommutative torus  ${\cal A}_{\theta}$  with  the 
unit ${1\over\mu}e$    is a coordinate ring of elliptic curve ${\cal E}_{\tau}$.   
Moreover,   such a correspondence is a covariant functor  which maps isomorphic 
elliptic curves to the   stably isomorphic (Morita equivalent) noncommutative tori; 
the latter fact follows from an observation that isomorphisms in category {\bf Mod} 
correspond to  stable isomorphisms in the category of underlying algebras. 
Such a  functor explains  the same  (modular)  transformation law in formulas  
(\ref{eq3}) and (\ref{eq5}).  
}
\end{rmk}
\begin{lem}\label{lem3}
The coordinate ring of elliptic curve ${\cal E}_{CM}^{(-d,f)}$ is isomorphic to 
the noncommutative torus ${\cal A}_{RM}^{(d, {\goth f})}$ with the unit ${1\over \log\varepsilon}e$,
where ${\goth f}$ is the least integer satisfying equation $|Cl~({\goth R}_{\goth f})|=|Cl~(R_f)|$
and $\varepsilon$ is the fundamental unit of order ${\goth R}_{\goth f}$. 
\end{lem}
{\it Proof.} 
The fact that ${\cal A}_{RM}^{(d, {\goth f})}$ is a coordinate ring of elliptic curve 
 ${\cal E}_{CM}^{(-d,f)}$ was proved in  \cite{Nik1}.  
We shall focus on the second part of lemma \ref{lem3}   saying that the scaling constant
 $\mu=\log\varepsilon$.  
To express $\mu$ in terms of intrinsic invariants of pseudo-lattice 
$K_0^+({\cal A}_{RM}^{(d, {\goth f})})\cong {\Bbb Z}+{\Bbb Z}\theta$,
recall that   ${\goth R}_{\goth f}$ is the  ring of endomorphisms of    
${\Bbb Z}+{\Bbb Z}\theta$;  we shall write ${\goth R}_{\goth f}^{\times}$ to denote
the multiplicative group of  units (i.e. invertibe elements) of ${\goth R}_{\goth f}$.  
Since $\mu$ is an additive functional on the  pseudo-lattice $\Lambda={\Bbb Z}+{\Bbb Z}\theta$,
for each   $\varepsilon, \varepsilon'\in {\goth R}_{\goth f}^{\times} $ it  must hold
$\mu(\varepsilon\varepsilon' \Lambda)=\mu(\varepsilon\varepsilon') \Lambda=
\mu(\varepsilon)\Lambda+\mu(\varepsilon')\Lambda$.
Eliminating  $\Lambda$ in the last equation,  one gets
\begin{equation}
\mu(\varepsilon\varepsilon')=\mu(\varepsilon)+\mu(\varepsilon'),
\qquad \forall \varepsilon, \varepsilon'\in    {\goth R}_{\goth f}^{\times}.
\end{equation}
The only real-valued function on ${\goth R}_{\goth f}^{\times}$ with such a property
is the logarithmic function (a regulator of ${\goth R}_{\goth f}^{\times}$);  thus $\mu(\varepsilon)=\log\varepsilon$,
where $\varepsilon$ is the fundamental unit of  ${\goth R}_{\goth f}$.   
Lemma \ref{lem3} is proved.
 $\square$

\begin{rmk}
{\bf (Second proof of lemma \ref{lem3})}
\textnormal{
The formula $\mu=\log\varepsilon$ can be derived  using  a purely measure-theoretic argument.  
Indeed,  if $h_x: {\Bbb R}\to {\Bbb R}$ is a ``stretch-out'' automorphism 
of real line ${\Bbb R}$ given by the formula $t\mapsto tx,~\forall t\in {\Bbb R}$,
then the only $h_x$-invariant measure $\mu$ on ${\Bbb R}$ is the ``scale-back''
measure $d\mu={1\over t} dt$.  Taking the antiderivative and integrating 
between $t_0=1$ and $t_1=x$,  one gets
\begin{equation}
\mu=\log x.
\end{equation}
It remains to notice that for pseudo-lattice 
$K_0^+({\cal A}_{RM}^{(d,{\goth f})})\cong {\Bbb Z}+{\Bbb Z}\theta$,
the automorphism $h_x$ corresponds to $x=\varepsilon$,  where $\varepsilon>1$
is the fundamental unit of  order ${\goth R}_{\goth f}$.  
Lemma \ref{lem3} follows. $\square$. 
}
\end{rmk}

\bigskip
One  can prove   theorem \ref{thm1} in the following steps.

\bigskip
(i)  Let  $d\not\in\{1,2,3,7,11,19,43, 67,163\}$ be a positive square-free integer.
In this case $h=|Cl~(R_f)|\ge 2$ and   ${\cal E}_{CM}^{(-d,f)}\not\cong {\cal E}({\Bbb Q})$.

\bigskip
(ii) Let $\{{\cal E}_1,\dots, {\cal E}_h\}$ be pairwise  non-isomorphic elliptic curves 
having the same endomorphism ring $R_f$.   From  $|Cl~(R_f)|=|Cl~({\goth R}_{\goth f})|$ and lemma \ref{lem3},  
one gets  $\{{\cal A}_1,\dots, {\cal A}_h\}$ pairwise stably non-isomorphic  noncommutative tori;  the corresponding
pseudo-lattices $K_0^+({\cal A}_i)={\Bbb Z}+{\Bbb Z}\theta_i$  will have  the same
endomorphism ring ${\goth R}_{\goth f}$.   Thus    for each $1\le i\le h$ one gets 
an inclusion 
\begin{equation}
(\log\varepsilon) e^{2\pi i\theta_i}\in H(k),
\end{equation}
where $H(k)$ is the Hilbert class field of quadratic field $k={\Bbb Q}(\sqrt{-d})$ modulo conductor $f$.  
Since $(\log\varepsilon)\exp (2\pi i\theta_i)=\exp (2\pi i\theta_i+\log\log\varepsilon):={\cal J}(\theta_i, \varepsilon)$,
one concludes that  ${\cal J}(\theta_i, \varepsilon)\in H(k)$.

\bigskip
(iii)   The  numbers ${\cal J}(\theta_i, \varepsilon)$ are algebraically conjugate.  Indeed,
the ideal class group $Cl~({\goth R}_{\goth f})$ acts transitively on the set of 
numbers $\theta_i$;  so will be its action on the set  ${\cal J}(\theta_i, \varepsilon)$. 
But  $Cl~({\goth R}_{\goth f})\cong Cl~(R_f)\cong Gal~(H(k)|k)$ and,  therefore,
the Galois group $Gal~(H(k)|k)$ acts transitively in the set of pairwise distinct algebraic
numbers ${\cal J}(\theta_i, \varepsilon)$.   The latter happens if and only if the  numbers
are algebraically conjugate.

\bigskip
Theorem \ref{thm1} is proved.
$\square$

\vskip1cm

\textsc{Department of Mathematics and Computer Science, St.~John's University, 8000 Utopia Parkway,  
New York,  NY 11439, United States;} ~\textsc{E-mail:} {\sf igor.v.nikolaev@gmail.com}

\end{document}